\documentclass[11pt, english]{article}

\usepackage[margin= 2.5 cm]{geometry}

\usepackage{amsthm}
\usepackage{amsmath}
\usepackage{amssymb}
\usepackage{setspace}
\usepackage{mathtools}
\usepackage{verbatim}
\usepackage{booktabs}
\usepackage{graphicx}
\usepackage[hidelinks]{hyperref}
\usepackage{cleveref}
\usepackage{graphicx}
\usepackage{appendix}
\usepackage[inline]{enumitem}
\usepackage{framed}
\usepackage{subcaption}
\usepackage{microtype,xurl}
\usepackage{caption}
\usepackage{listings}

\usepackage{floatrow}
\usepackage[T1]{fontenc}

\lstdefinestyle{verifier}{
  language=Python,
  basicstyle=\ttfamily\fontsize{8.1}{10.2}\selectfont,
  keywordstyle=\bfseries,
  commentstyle=\itshape,
  stringstyle=\ttfamily,
  showstringspaces=false,
  keepspaces=true,
  columns=fullflexible,
  breaklines=true,
  breakatwhitespace=false,
  numbers=left,
  numberstyle=\tiny,
  numbersep=7pt,
  xleftmargin=1.6em,
  aboveskip=0.8em,
  belowskip=0.8em,
  tabsize=4,
  upquote=true
}

\theoremstyle{plain}

\newtheorem*{thm*}{Theorem}
\newtheorem{thm}{Theorem}
\Crefname{thm}{Theorem}{Theorems}
\newtheorem*{lem*}{Lemma}

\newtheorem*{claim*}{Claim}

\Crefname{claim}{Claim}{Claims}
\Crefname{claim}{Claim}{Claims}

\Crefname{prop}{Proposition}{Propositions}

\Crefname{cor}{Corollary}{Corollaries}

\Crefname{conj}{Conjecture}{Conjectures}

\Crefname{qn}{Question}{Questions}

\Crefname{obs}{Observation}{Observations}

\newtheorem{ex}[thm]{Example}
\Crefname{ex}{Example}{Examples}

\theoremstyle{definition}

\Crefname{prob}{Problem}{Problems}

\Crefname{defn}{Definition}{Definitions}

\newtheorem*{defn*}{Definition}

\theoremstyle{remark}

\expandafter\def\expandafter\normalsize\expandafter{%
    \normalsize
    \setlength\abovedisplayskip{8pt}
    \setlength\belowdisplayskip{8pt}
    \setlength\abovedisplayshortskip{4pt}
    \setlength\belowdisplayshortskip{4pt}
}

\usepackage[square,sort,comma,numbers]{natbib}
 \setlist[itemize]{leftmargin=*}

 \newcommand{\EE}{\mathbb E}
 \newcommand{\R}{\mathbb R}
 \newcommand{\one}{\mathbf1}
 \newcommand{\file}[1]{\mbox{\nolinkurl{#1}}}
\renewcommand{\ex}{\operatorname{ex}}

\begin{document}
\title{The uniform Tur\'an density of the tetrahedron}
\author{
Matija Buci\'c\thanks{Faculty of Mathematics, University of Vienna, Vienna, Austria. Email: \href{mailto:matija.bucic@univie.ac.at} \textbf{matija.bucic@univie.ac.at}. }
% \and
% Haoran Luo\thanks{Department of Mathematics, Statistics and Computer Science, University of Illinois Chicago, Chicago, Illinois, USA. Email: \textbf{haoranl8@uic.edu}. The author was partially supported by an AMS-Simons Travel Grant.}
}
\date{}
\maketitle
\begin{abstract}
We prove that the uniform Tur\'an density of $K_4^{(3)}$ is equal to $1/2$.  Kielak, Král', Lamaison, Liu, Shu, and Wu have recently proved the same result using combinatorial methods, while our proof is Fourier-analytic. This result answers a question of Erd\H{o}s and S\'os from the founding 1982 paper on this topic.
\end{abstract}

\section{The combinatorial statement}

A central problem in extremal combinatorics is to determine how many edges
a graph or hypergraph can have without containing a prescribed subgraph
or subhypergraph. For a $k$-uniform hypergraph $F$, the \emph{Tur\'an number}
$\ex(n,F)$ is the maximum number of edges in an $n$-vertex $k$-uniform
hypergraph containing no copy of $F$. The asymptotic behaviour of this
quantity is described by the \emph{Tur\'an density} $\pi(F):=\lim_{n\to\infty}\frac{\ex(n,F)}{\binom{n}{k}},$
whose existence follows from a classical result of
Katona, Nemetz, and Simonovits~\cite{KatNS64}. Thus, $\pi(F)$ is the largest
edge density asymptotically attainable by $F$-free hypergraphs.
Equivalently, it is the infimum of all $d$ such that every sufficiently
large $k$-uniform hypergraph of edge density at least $d$ contains $F$.

In the special case of graphs, these questions are substantially better understood.
The theorems of Mantel~\cite{Man07} and Tur\'an~\cite{Tur41} determine
$\ex(n,F)$ exactly when $F$ is a complete graph, while the
Erd\H{o}s--Stone theorem~\cite{ErdS46} determines the Tur\'an density
of every graph $F$ with chromatic number $r\ge 2$ (see also~\cite{ErdS66}).
For hypergraphs of uniformity at least three,
however, the corresponding problems are considerably more difficult.
Erd\H{o}s~\cite{Erd81} offered \$1,000 for determining the Tur\'an densities of all complete $k$-uniform hypergraphs with $k\ge 3$, and \$500 for
determining the density of even one such hypergraph on at least $k+1$
vertices. The difficulty is already evident in the smallest nontrivial
case, the complete $3$-uniform hypergraph $K_4^{(3)}$ on four vertices, usually referred to as the \emph{tetrahedron}.
Determining its Tur\'an density has been a central problem since
Tur\'an's work in 1941~\cite{Tur41}, motivating substantial research,
including~\cite{FraF84,ChuL99} and bounds obtained using Razborov's
flag algebra method~\cite{BabT11,Raz10,Raz07}. For broader accounts
of hypergraph Tur\'an problems, we refer the reader to the surveys
of Sidorenko~\cite{Sid95}, Keevash~\cite{Kee11}, and Balogh, Clemen, and Lidick\'{y}~\cite{BalCL21}.

Tur\'an conjectured that $\pi(K_4^{(3)})=5/9$. The following construction
establishes the lower bound. Partition an $n$-vertex set into three
parts $V_1,V_2,V_3$ of sizes as equal as possible, and include every
triple that either meets all three parts or has two vertices in $V_i$
and one vertex in $V_{i+1}$ for some $i\in\{1,2,3\}$, with subscripts
taken modulo $3$. The resulting hypergraph is $K_4^{(3)}$-free and
has edge density tending to $5/9$ as $n\to\infty$. Its edges are,
however, distributed very unevenly: each part $V_i$ is an independent
set, so a positive global edge density coexists with large vertex
subsets spanning no edges.

This phenomenon motivated Erd\H{o}s and S\'os~\cite{ErdS82,Erd90}
to study Tur\'an problems under an additional requirement that edges
be distributed uniformly across large vertex subsets. A three-uniform hypergraph $G$ on $N$ vertices is \emph{uniformly
$(d,\eta)$-dense} if
\[
 e_G(U)\ge d\binom{|U|}{3}-\eta N^3
 \qquad\text{for every }U\subseteq V(G).
\]
The \emph{uniform Tur\'an density} $\pi_{\mathrm u}(F)$ is the supremum of
those $d\in[0,1]$ for which, for every $\eta>0$, there are arbitrarily large
$F$-free uniformly $(d,\eta)$-dense three-uniform hypergraphs.

Unlike ordinary Tur\'an
density, which imposes only a global density requirement, uniform
Tur\'an density requires the same lower bound on the edge density
of every linear-sized vertex subset.

Among the problems posed by Erd\H{o}s and S\'os in their foundational 1982 paper~\cite{ErdS82} were the determination
of the uniform Tur\'an densities of $K_4^{(3)}$ and $K_4^{(3)-}$,
where $K_4^{(3)-}$ is obtained from $K_4^{(3)}$ by deleting one edge.
These examples illustrate the difficulty of uniform Tur\'an problems
even for small $3$-uniform hypergraphs. The equality
\[
    \pi_{\mathrm{u}}(K_4^{(3)-})=\frac14
\]
was established in~\cite{GleKV16,ReiRS18a}, and since then, numerous papers have explored the topic with many great success stories in recent years ~\cite{ReiRS18,BucCKMM23,GarKL21,lamaison,LamW24,KraKLT25,LiLWZ25,DinLLWY25,LinZ25,LinWZ25,CheS26,GarIKKL26,LinSWZ26,LiLP26}. The corresponding problem
for $K_4^{(3)}$ has proved substantially more difficult, with essentially every paper on the topic since the foundational one by Erd\H{o}s and S\'os mentioning it as one of the central open problems in the area. In parallel with Kielak, Král', Lamaison, Liu, Shu, and Wu \cite{dan-paper}, we developed a different approach that yields an alternative solution to this problem.

\begin{thm}\label{thm:main}
For the complete three-uniform hypergraph on four vertices $K_4^{(3)}$, we have
\[
 \pi_{\mathrm u}(K_4^{(3)})=\frac12.
\]
\end{thm}

We note that the lower bound comes from a very natural construction of R\"odl from 1986~\cite{Rod86}, a variant of which we will present in the following section.

\section{The palette reduction}

In the proof, we will implicitly make use of the hypergraph regularity setup introduced by Reiher, R\"odl, and Schacht \cite{ReiRS18a} combined with a remarkable recent result of Lamaison \cite{lamaison}, which reduces the general question of determining the Turán density of an arbitrary fixed $3$-uniform hypergraph to only verifying so-called palette constructions. For readers unfamiliar with the topic, we note that this is similar to the graph regularity setup, where one only needs to solve the problem in the (finite) reduced graph. We will not introduce the setup in full generality here, and rather point a reader new to the topic to the survey \cite{Rei20} and the paper of Lamaison itself \cite{lamaison}.

\begin{defn*}
A \emph{palette} is a set $P\subseteq C^3$, where $C$ is a finite nonempty
set, whose elements are called colours.  Its \emph{density} is $|P|/|C|^3$.
We say that $P$ \emph{admits a tetrahedron} if there are colours
$a,b,c,e,f,g\in C$ for which
\begin{equation}\label{eq:configuration}
 (a,b,e),\quad(a,c,f),\quad(b,c,g),\quad(e,f,g)\ \in P.
\end{equation}
A palette admitting no such assignment is said to be
\emph{tetrahedron-free}.
\end{defn*}

We prove the following upper bound on the density of any tetrahedron-free palette, which, combined with Lamaison's palette theorem~\cite[Theorem~1.1]{lamaison}, specialized to the case of the tetrahedron, immediately implies \Cref{thm:main}.

\begin{thm}\label{lem:palette}
Every tetrahedron-free palette $P\subseteq C^3$ satisfies $|P|\le\frac12|C|^3$.
\end{thm}

The well-known construction $P=\{(x,y,z)\in \{0,1\}^3:x\ne y\}$ mentioned above establishes tightness.

\section{Proof of the palette bound}\label{sec:proof}
In this section, we prove Theorem~\ref{lem:palette}.
Let $P\subseteq C^3$ be tetrahedron-free, put $n=|C|$, and write
\[
 p:=\one_P,\qquad d:=\frac{|P|}{n^3},
\]
where $\one_P:C^3 \to \{0,1\}$ denotes the indicator function of $P$.

All expectations below use independent uniform colours from $C$.
The normalized number of tetrahedra in our pallete (recall~\eqref{eq:configuration}) is given by
\begin{equation}\label{eq:count}
 t(P)=\EE_{a,b,c,e,f,g}\,
 p(a,b,e)p(a,c,f)p(b,c,g)p(e,f,g) %=0.
\end{equation}

The case $n=|C|=1$ of \Cref{lem:palette} is immediate: the only pattern, if present, would itself give a tetrahedron.  Hence, we assume $n\ge2$.% and set $m=n-1$.

\subsection{Expanding the indicator}
Choose a real orthonormal basis $\phi_0,\ldots,\phi_{n-1}$ for functions from $C \to \R$
with respect to the normalized inner product $\langle f, g \rangle:= \EE_x f(x)g(x) = \frac1n \sum_{x\in C} f(x)g(x)$, with $\phi_0\equiv1$.
Thus,
\[
 \EE_x\phi_i(x)\phi_j(x)=\one_{\{i=j\}}.
\]
Write the orthonormal expansion of the palette indicator function $p:C^3 \to \mathbb R$ as
\begin{equation}\label{eq:expansion}
 p(x,y,z)=\sum_{i,j,k=0}^{n-1}h_{ijk}\phi_i(x)\phi_j(y)\phi_k(z),
\quad \textrm{where} \,\,\,
h_{ijk}=\EE_{x,y,z}p(x,y,z)\phi_i(x)\phi_j(y)\phi_k(z).
\end{equation}
Note that since $\phi_0 \equiv 1$, we have $h_{000}=d$.  Parseval's identity and $p^2=p$ (being an indicator)
give
\begin{equation} \label{eq:variance}
 \sum_{i,j,k=0}^{n-1}h_{ijk}^2=\EE p^2=d,
 \qquad
 V:=\sum_{(i,j,k)\ne(0,0,0)}h_{ijk}^2=d-d^2.
\end{equation}
Next, let us substitute~\eqref{eq:expansion} into~\eqref{eq:count}.  Each colour variable occurs in exactly two factors. Changing the order of summation and using orthogonality gives
\begin{equation}\label{eq:coefficientcount}
 t(P)=\sum_{a,b,c,e,f,g=0}^{n-1}
             h_{abe}h_{acf}h_{bcg}h_{efg}.
\end{equation}

\subsection{The finite positivity certificate}
The final idea is to, in some sense, ``complete'' the square where the expressions in \eqref{eq:variance} give us the control of the sum of squares of $h_{ijk}$'s and expressions in \eqref{eq:coefficientcount} give the ``cross-terms''.

In particular, the certificate described in Appendix~\ref{app:certificate} gives the exact identity
\begin{equation}\label{eq:maincertificate}
 10^{12}\bigl(4t(P)+3V^2-3d^4\bigr)
 =\sum_{r=0}^{6}\ \sum_{\boldsymbol i\in\{1,\ldots,n-1\}^{r}}
       q_r(\boldsymbol i)^{\mathsf T}R_rq_r(\boldsymbol i),
\end{equation}
where the entries of $q_r(\boldsymbol i)$ are specified quadratic
expressions in the coefficients $h_{ijk}$, and the seven fixed integer
matrices $R_r$ (of finite orders, each at most $360$) are positive definite.
%For $r=0$ the inner sum has one term.
The identity and positive definiteness are checked exactly by the
supplied verifier; Appendix~\ref{app:certificate} explains the verification process.

By positive definiteness, every term on the right of~\eqref{eq:maincertificate} is nonnegative.
Since $P$ is tetrahedron-free, we have $t(P)=0$. Also by $V=d-d^2$, it follows that
\[
 0\le 3(d-d^2)^2-3d^4=3d^2(1-2d).
\]
This implies $d\le1/2$, as desired.  This proves Lemma~\ref{lem:palette} and hence
Theorem~\ref{thm:main}.
\qed

\section{Concluding remarks}\label{sec:conc-remarks}

In this paper, we determine the uniform Tur\'an density of the tetrahedron, arguably one of the most well-known and classical open problems in uniform Tur\'an theory.

We note that \Cref{thm:main} in fact implies that uniform Tur\'an density of an infinite family of graphs is equal to a half by a simple lifting argument based on Lamaison's palette theorem. Specifically, any $3$ uniform hypergraph whose links are not bipartite (forcing the uniform Tur\'an density to be at least $1/2$ using the standard construction) and which admits a $4$-coloring of its vertices in which every hyperedge receives 3 different colors has uniform Tur\'an density equal to a half. Basically, the latter condition allows us to embed our hypergraph into a ``blow-up'' of the tetrahedron which the main theorem, combined with the pallete theorem implies exists. Another notable example that falls into this category is the Icosahedron hypergraph with edges being triples of vertices defining its faces.
More generally, for an odd integer $\ell \ge 5$ let $I_\ell$ be the hypergraph with vertex set $\{b,v_1,\ldots, v_{\ell},u_1,\ldots, u_{\ell},t\}$ and edge set consisting of all triples of the form $\{v,v_i,v_{i+1}\}$, $\{t,u_i,u_{i+1}\}$, $\{v_i,u_i, v_{i+1}\}$, and $\{v_i,u_i, u_{i-1}\}$, where the indices are taken modulo $\ell$. Here, $I_5$ is precisely the icosahedron hypergraph. It is easy to see that this family has non-bipartite links and that it admits a $4$-coloring, so $\pi_u(I_\ell)=1/2$ for any $\ell \ge 5$.

The approach of proving uniform Tur\'an density bounds using Fourier inspired ideas we bring to the table has a great potential for determining the Tur\'an density of other hypergraphs. In some sense our argument here is a second moment one, which requires that every pair of vertices belongs to at most two hyperedges. This is a property which fails for larger cliques, namely to attack $K_5^{(3)}$ one would need in a certain sense a third moment argument and the complexity of that seems to go through the roof. That said there are other families where a similar argument goes through. For example, the simple proof of $\pi_u(K_4^{(3)-})=1/4$ from \cite[Section 6.1]{lamaison} can be rephrased in terms of establishing an equality along the lines of \eqref{eq:maincertificate}, but requiring a much simpler and easily hand verified RHS. This reformulation was behind one of our original ideas of how to attack the full tetrahedron case. The same strategy can also recover the main result from \cite{GarIKKL26} showing that a certain infinite family of hypergraphs has uniform Tur\'an density equal to $8/27$, using a similar, simple, easily hand verifiable equality.

The computer-assisted portion of our proof leaves a lot to be desired, in particular given that it obscures the actual mathematical reasons behind the result.
For these reasons, it would be of great interest to find a simple proof not requiring computer assistance.

\textbf{Acknowledgments.} The first author would like to thank Samuel Mohr, Nina Kam\v{c}ev, and Haoran Luo for useful discussions.

\textbf{Declaration of use of AI}. The identity \eqref{eq:maincertificate} was found by ChatGPT 6 Pro, following an upload of an early draft containing various ideas the author had for attacking the problem, and following a prolonged discussion in which the author guided its progress.
Eventually, ChatGPT proposed a more general tensor orthogonal projection inequality (hiding behind \eqref{eq:maincertificate}) and wrote a number of computer programs to find and generate the expansion as a sum of positive definite quadratic forms in order to prove the desired inequality. The code that led to the inequality and the guiding principles for finding it are available upon request, but since they are technically not required for the proof, once the identity \eqref{eq:maincertificate} is identified, we do not include them here. Besides the author hand-checking the correctness of the verifier, Aristotle was used to verify its soundness both in Lean and to run a separate audit. The equality \eqref{eq:maincertificate} was too large to directly verify in Lean.

\textbf{Note.}
The same result was obtained in \cite{dan-paper}. The second author of \cite{dan-paper} shared with us on August 18, 2026, that he and his collaborators had obtained a computer-assisted proof of the result that does not contain AI-assisted arguments. We disclosed that we had also been working independently on the problem and expressed our already existing intention to feed our strategies into an LLM. ChatGPT-5.6, which we subsequently used, was unable to complete the proof of the result. On September 6, 2026, two days after its release, ChatGPT-6 managed to complete the proof, using an additional prompt added to the previous conversation with ChatGPT-5.6. On September 8, 2026, after verifying and simplifying the argument, we informed the authors of \cite{dan-paper} about having obtained an alternative computer-assisted proof and shared a draft - at that point, their manuscript had already been submitted to arXiv, but had not yet become publicly available or been shared with us. The two proofs use different approaches.

\providecommand{\MR}[1]{}
\providecommand{\MRhref}[2]{%
  \href{http://www.ams.org/mathscinet-getitem?mr=#1}{#2}
}

   \bibliographystyle{amsplain_initials_nobysame}
   \bibliography{ref}

\appendix
\section{Certificate and exact verification}\label{app:certificate}

This appendix specifies the objects in~\eqref{eq:maincertificate} and the
finite checks that establish it.  The two proof-bearing files
are \file{certificate.json} and \file{verify.py}.  The former contains
integer data; the latter checks a formal identity and seven matrix
positivity certificates.  The program's variable $T$ denotes the
coefficient array $h$ from~\eqref{eq:expansion}.  Throughout this appendix,
put $m=n-1$, so that every coefficient index is either $0$ or an element
of $\{1,\ldots,m\}$.

The certificate is an explicit way of writing the expression we wish to
prove nonnegative as a sum of nonnegative quadratic forms.  To check such
a representation, we need to answer two questions:
\begin{enumerate}
\item Does expanding the proposed representation give exactly
$10^{12}(4t(P)+3V^2-3d^4)$?
\item Is each quadratic form in that representation nonnegative for
every real choice of its arguments?
\end{enumerate}
The first question is a matter of expanding products and collecting
like terms.  The second is a finite matrix calculation, followed by the
elementary diagonal-dominance argument given below.  The data file
supplies the expressions and matrices needed for these calculations;
the verifier performs the calculations afresh.  In particular, finding
the certificate and checking it are separate tasks: once the data are
given, checking them requires only the operations described here.

For these checks, the coefficients $h_{ijk}$ can be regarded as arbitrary
real variables.  The program works with formal expressions in these
variables, keeping the sums over indices symbolic.  The facts that
$h_{000}=d$, that $V=d-d^2$ for an indicator, and that $t(P)=0$ for a
tetrahedron-free palette are used when applying the resulting inequality
in Section~\ref{sec:proof}.  The factor $10^{12}$ is an exact scaling
factor which allows the identity to be checked using integer
coefficients; it specifies no numerical tolerance.

\subsection{Reading the data}
An entry of $q_r(\boldsymbol i)$ is quadratic in the coefficients $h_{ijk}$,
so it is built from products of two coefficients, possibly with some
indices summed out.  The number $r$ records how many indices are left
free at this stage.  There are at most six such indices, because two
coefficients have six index positions in total.  This explains the seven
blocks $r=0,\ldots,6$ in~\eqref{eq:maincertificate}.  For a fixed tuple
$\boldsymbol i$, the vector $q_r(\boldsymbol i)$ has $b_r$ real entries;
the block sums $q_r(\boldsymbol i)^{\mathsf T}R_rq_r(\boldsymbol i)$ over
all $m^r$ choices of that tuple.  When $r=0$,
there is one empty tuple, so the block contributes one quadratic form.
Since the entries of $q_r$ are themselves quadratic in $h$, a quadratic
form in these entries has degree four in $h$, matching the expression
on the left of~\eqref{eq:maincertificate}.

For each $r=0,\ldots,6$, a \emph{raw expression} is encoded by two ordered
triples of nonnegative integer labels.  Label $0$ denotes the fixed
coefficient index $0$.  The labels $1,\ldots,r$ occur once each and denote
the free indices $i_1,\ldots,i_r$.  Each larger label occurs exactly twice
and denotes an index summed independently over $1,\ldots,m$.  The
expression is the product of the two indicated coefficients, summed over
all these internal indices.  For example, with $r=2$,
\[
 \begin{aligned}
 ((0,1,2),(0,0,0))&\quad\longmapsto\quad h_{0ij}h_{000},\\
 ((1,3,0),(2,3,0))&\quad\longmapsto\quad
                    \sum_{k=1}^{m}h_{ik0}h_{jk0}.
 \end{aligned}
\]
All positive labels are contiguous.  No distinctness is imposed on the
values of any indices.  The labels are names of indices, rather than
their numerical values.  Thus, in the second example, labels $1$ and $2$
stand for the freely chosen values $i$ and $j$, while the two occurrences
of label $3$ instruct us to use the same summation index $k$ in both
factors.  The cases $i=j$, $i=k$, or $j=k$ are all included.

Block $r$ lists its raw expressions $F_{r,0},\ldots,F_{r,k_r-1}$ in
\file{features}.  Their integer linear combinations are recorded in
\file{combination_columns} and define the $b_r$ entries of $q_r$:
\begin{equation}\label{eq:q}
 (q_r(\boldsymbol i))_u
 =\sum_{j=0}^{k_r-1}N^{(r)}_{ju}F_{r,j}(\boldsymbol i),
 \qquad 0\le u<b_r.
\end{equation}
Each column of $N^{(r)}$ is stored sparsely as pairs $[j,c]$, meaning
coefficient $c$ times raw expression $j$.  For example, a column
$[[2,3],[5,-1]]$ would specify the entry $3F_{r,2}-F_{r,5}$; coefficients
of unlisted raw expressions are zero.  Thus \file{features} gives the
building blocks, and \file{combination_columns} says how to assemble
them into the entries of $q_r$.

The field \file{gram_lower}
stores the lower triangle of $R_r$, completed by symmetry.  The field
\file{congruence_lower} stores a lower-triangular integer matrix $B_r$,
completed by zeros above the diagonal.  Both matrices have order $b_r$.
The matrix $R_r$ supplies the coefficients of the quadratic form
$q_r^{\mathsf T}R_rq_r$.  The auxiliary matrix $B_r$ is used only to prove
that this form is positive definite; it does not enter the expansion
used to check the identity.

\begin{center}
\begin{tabular}{@{}lrrrrrrr@{}}
\toprule
$r$ & 0 & 1 & 2 & 3 & 4 & 5 & 6\\
\midrule
$k_r$ & 24 & 39 & 87 & 168 & 288 & 360 & 360\\
$b_r$ & 22 & 39 & 85 & 168 & 286 & 360 & 360\\
\bottomrule
\end{tabular}
\end{center}

The file's format identifier is \file{tetrahedron-tensor-sos-v1} and its
\file{denominator} is $10^{12}$.  The verifier checks dimensions,
integrality, and the stated label conventions as it reads each block.
The numbers $k_r$ and $b_r$ describe the sizes of the supplied data;
they do not depend on $m$.

\subsection{Checking the identity}
To check an ordinary polynomial identity, one expands both sides and
compares the coefficients of each monomial.  Here the same idea applies,
except that a single term may stand for a sum of many monomials over
indices ranging from $1$ to $m$.  The verifier records the pattern of
these sums, so it can compare the two sides without choosing a value
of $m$.

After expanding~\eqref{eq:q}, each term on the right
of~\eqref{eq:maincertificate} is a sum of products of four coefficients
$h_{ijk}$.  It is encoded by four ordered triples of labels.  Each positive
label occurs twice and is summed over $1,\ldots,m$; label $0$ is fixed.
An internal index already occurs twice in its own raw expression, while
a free index occurs once in each of the two raw expressions being
multiplied.  Thus, after summing over the free indices, every positive
label occurs exactly twice.
Two operations preserve the value of such an expression: permuting the
four scalar factors, and renaming dummy summation indices.  Neither
operation permutes the three positions within a coefficient $h_{ijk}$.
This restriction matters because the palette need not be symmetric:
for example, $h_{ijk}$ and $h_{jik}$ need not be equal.

To see the index bookkeeping in a concrete case, write
$G(i,j)=\sum_{k=1}^{m}h_{ik0}h_{jk0}$ for the second raw expression above.
Its contribution when paired with itself is
\[
 \sum_{i,j=1}^{m}G(i,j)^2
 =\sum_{i,j,k,\ell=1}^{m}
       h_{ik0}h_{jk0}h_{i\ell0}h_{j\ell0}.
\]
The two copies share the free indices $i,j$, but have separate internal
summation indices $k,\ell$.  Using the same internal label in both
copies would incorrectly keep only the terms with $k=\ell$.  Separate
labels mean that the two indices are summed independently; their values
are still allowed to coincide.  Before canonicalization, the four
factors in this example are encoded by
\[
 ((1,3,0),(2,3,0),(1,4,0),(2,4,0)).
\]

The function \file{canonical_graph} takes the lexicographically least
encoding under these operations, considering all $4!$ factor permutations
and renaming positive labels in order of first occurrence.  In other
words, it gives each sum a standard written form, so that changing the
names of its dummy indices or the order of its factors does not produce
a different entry in the coefficient list.  The name \file{canonical_graph}
reflects the pattern of pairings between index positions: each positive
label joins its two occurrences.  The function
\file{inner_product_graph} expands a product of two raw expressions with
the same free indices.  It identifies their free labels and keeps their
internal labels disjoint.  These are precisely the index identifications
arising when summing over $\boldsymbol i$ in~\eqref{eq:maincertificate}.
The verifier accumulates the resulting coefficients as integers.  Since
$R_r$ is symmetric, it counts diagonal matrix entries once and
above-diagonal entries twice, just as in the expansion
\[
 q_r^{\mathsf T}R_rq_r
 =\sum_u (R_r)_{uu}(q_r)_u^2
   +2\sum_{u<v}(R_r)_{uv}(q_r)_u(q_r)_v.
\]
Substituting~\eqref{eq:q} into this expression explains the nested loops
over entries of $R_r$ and columns of $N^{(r)}$ in \file{verify}.  Each
product of raw expressions contributes its integer coefficient to the
entry selected by \file{inner_product_graph}.

The left side is generated independently of the data.  For the count
in~\eqref{eq:coefficientcount}, begin with the four triples
\[
 (1,2,4),\quad(1,3,5),\quad(2,3,6),\quad(4,5,6).
\]
Expand the $2^6$ possibilities for each label to be either $0$ or a
positive summation index.  For $V^2$, group the coefficients in $V$ by the
seven nonempty subsets of their three positions that have positive
indices, and expand the $7^2$ ordered products of their squared-norm
sums.  Use disjoint dummy indices in the two factors.  Finally subtract
$3h_{000}^4$, with the factors $4$, $3$, and $10^{12}$ as
in~\eqref{eq:maincertificate}.

For example, the seven groups contributing to $V$ include
$\sum_i h_{i00}^2$ and $\sum_{i,j}h_{ij0}^2$, as well as the analogous
sums for the other choices of positive positions.  The all-zero choice
is excluded because $V$ omits $h_{000}^2$.  This is why the code uses
seven choices for each factor of $V^2$, while the tetrahedron count has
two choices for each of its six indices.  The function
\file{target_polynomial} builds this expansion directly from the formulas
in the main proof.  The coefficient lists called \file{actual} and
\file{target} therefore come from the proposed representation and the
expression it is supposed to represent, respectively.

All $2\,752$ coefficients in the union of the two resulting lists agree
exactly.  This proves~\eqref{eq:maincertificate}: every identification used
in the comparison preserves the corresponding finite sum.  Linear
independence of the canonical expressions is unnecessary.  The procedure
never fixes $m$, so the same identity holds for every number of colours.
To spell out the last point, after putting both sides into the same
notation, their difference has coefficient zero in front of every
recorded sum.  It is therefore zero for every array $h$, regardless of
whether additional relations between those sums hold in a particular
dimension.  A sum present on only one side is also checked, with its
coefficient on the other side taken to be zero.  In particular, all
additional terms introduced by expanding the certificate must cancel
exactly.

\subsection{Checking positivity}
The identity check tells us what the quadratic forms add up to.  We
still need to know that each form is nonnegative.  The matrix $B_r$
supplied with each block provides a change of coordinates in which
positivity can be read from inequalities between integer entries.
This is the role of the second part of the certificate.

For each block the verifier checks that the diagonal of $B_r$ is positive
and computes, using integers,
\[
 H_r=B_rR_rB_r^{\mathsf T}.
\]
It checks symmetry and the strict diagonal-dominance inequalities
\begin{equation}\label{eq:dominance}
 (H_r)_{ii}>\sum_{j\ne i}|(H_r)_{ij}|.
\end{equation}
In words, each diagonal entry is larger than the sum of the absolute
values of all other entries in its row.  The positive square terms in
the quadratic form therefore have enough weight to absorb every
possibly negative cross-term, with a positive amount left over.
All seven blocks pass.  More precisely, for any real vector $x$, the inequality
$2|x_ix_j|\le x_i^2+x_j^2$ gives
\[
 x^{\mathsf T}H_rx
 \ge\sum_i\left((H_r)_{ii}-\sum_{j\ne i}|(H_r)_{ij}|\right)x_i^2.
\]
Thus $H_r$ is positive definite.  The triangular matrix $B_r$ is
invertible because its diagonal entries are nonzero.  For any nonzero
vector $y$, we can consequently write $y=B_r^{\mathsf T}x$ with $x\ne0$,
and then
\[
 y^{\mathsf T}R_ry
 =x^{\mathsf T}B_rR_rB_r^{\mathsf T}x
 =x^{\mathsf T}H_rx>0.
\]
Hence $R_r$ is positive definite as well.  This explains why checking
$H_r$ suffices.  The function \file{verify_positive_definite} carries out
the matrix multiplication and tests the row margins in~\eqref{eq:dominance}
using integers.  Taking $y=q_r(\boldsymbol i)$ now gives a nonnegative
term for every tuple $\boldsymbol i$, as required in Section~\ref{sec:proof}.

\subsection{Reproducing the finite check}
Place the data and verifier in the same directory and run
\begin{verbatim}
python3 -S verify.py certificate.json
\end{verbatim}
This command uses only the Python standard library.  All mathematical
checks use arbitrary-precision integers.  Running without \texttt{-S}
also permits an optional NumPy implementation of matrix multiplication
with Python-integer object arrays; no floating-point calculation enters
either version.  No optimization solver is required.

For a reader following the source in Appendix~\ref{app:verifier}, the
main routine \file{verify} processes the seven blocks in order.  In each
block it checks the data, proves positivity, and adds that block's
expansion to \file{actual}.  After all blocks have been processed, it
calls \file{target_polynomial} and compares every coefficient of
\file{actual} with \file{denominator} times the corresponding coefficient
of \file{target}.  The code raises an error if a required check fails.

A successful run prints one positivity confirmation for each block,
followed by confirmation that all $2\,752$ formal coefficients match
and the line beginning \file{VERIFIED}.  The printed expression $K(T)$
is the tetrahedral contraction in~\eqref{eq:coefficientcount}, and
$\|T\|_F^2-T_{000}^2$ is the sum of the squares of all coefficients except
$h_{000}$.  Thus, the final message is exactly the tensor inequality used
in the proof.  Its justification consists of the coefficient comparison
and the positivity checks above, together with the explicit interpretation
of labels as finite sums.

% For identification, the data file has $3\,782\,038$ bytes and the files have
% the following SHA-256 digests.
% \begin{center}
% \small
% \begin{tabular}{@{}ll@{}}
% \toprule
% File & SHA-256 digest\\
% \midrule
% \file{certificate.json} & \texttt{b5b651d0bfeef493f902f7b9d8dffdbbb746e2b41b97dc37215da9a553a07614}\\[0.3em]
% \file{verify.py} & \texttt{a570f18951947f8cb60887bb3f784c5bf40d9af774b1feee50a7924ef739fcdc}\\
% \bottomrule
% \end{tabular}
% \end{center}
%The complete verifier is supplied as executable source in the electronic supplement.  The data and program, not a saved output log, are the finite computational part of the proof.
%\clearpage

For identification, the data file has $3\,782\,038$ bytes and the files have
the following SHA-256 digests.
\begin{center}
\small
\begin{tabular}{@{}ll@{}}
\toprule
File & SHA-256 digest\\
\midrule
\file{certificate.json} & \texttt{b5b651d0bfeef493f902f7b9d8dffdbbb746e2b41b97dc37215da9a553a07614}\\[0.3em]
\file{verify.py} & \texttt{a570f18951947f8cb60887bb3f784c5bf40d9af774b1feee50a7924ef739fcdc}\\
\bottomrule
\end{tabular}
\end{center}
The complete verifier is supplied as executable source in the electronic supplement.  The data and program, not a saved output log, are the finite computational part of the proof.
\clearpage
\section{Verifier}\label{app:verifier}
%\noindent File: \file{experiments/three_projection_cut_certificates/verify.py}.
\begin{lstlisting}[style=verifier]
#!/usr/bin/env python3
"""Verify a dimension-independent tetrahedron tensor sum-of-squares identity.

Every mathematical check uses integers, not floating-point arithmetic. NumPy,
when available, only accelerates Python-integer matrix multiplication through
object arrays. There is a standard-library-only fallback.

Usage: python verify.py [certificate.json]
"""
from __future__ import annotations

import argparse
import itertools
import json
import time
from collections import Counter
from functools import lru_cache
from pathlib import Path
from typing import Sequence

try:
    import numpy as np
except ImportError:
    np = None

Row = tuple[int, int, int]
Feature = tuple[Row, Row]
Graph = tuple[Row, Row, Row, Row]


def fail_unless(condition: bool, message: str) -> None:
    if not condition:
        raise ValueError(message)


@lru_cache(maxsize=None)
def canonical_graph(rows: tuple[Row, ...]) -> tuple[Row, ...]:
    """Only commute the four tensor factors and rename dummy indices.

    Zero always means the distinguished coordinate; it is never renamed.
    The three slots within a tensor factor are NEVER permuted.
    """
    best = None
    for order in itertools.permutations(range(4)):
        renaming = {0: 0}
        result = []
        for factor in order:
            row = []
            for label in rows[factor]:
                if label not in renaming:
                    renaming[label] = len(renaming)
                row.append(renaming[label])
            result.append(tuple(row))
        candidate = tuple(result)
        if best is None or candidate < best:
            best = candidate
    return best


@lru_cache(maxsize=None)
def inner_product_graph(r: int, left: Feature, right: Feature) -> Graph:
    """Graph for the inner product of two quadratic rank-r tensor features.

    Free labels 1,...,r are shared. Contracted labels are made disjoint.
    """
    highest_left_label = max(max(row) for row in left)
    shift = highest_left_label - r
    right_renamed = tuple(
        tuple(label if label <= r else label + shift for label in row)
        for row in right
    )
    return canonical_graph(left + right_renamed)


def validate_feature(r: int, feature: Feature) -> None:
    fail_unless(len(feature) == 2 and all(len(row) == 3 for row in feature),
                "Every feature must contain two three-slot tensor factors.")
    labels = Counter(label for row in feature for label in row if label)
    fail_unless(all(isinstance(x, int) and x >= 0 for row in feature for x in row),
                "Feature labels must be nonnegative integers.")
    fail_unless(all(labels[i] == 1 for i in range(1, r + 1)),
                "Each free label must occur exactly once.")
    fail_unless(all(count == 2 for label, count in labels.items() if label > r),
                "Each contracted label must occur exactly twice.")
    fail_unless(set(labels) == set(range(1, max(labels, default=0) + 1)),
                "Feature labels must be contiguous.")


def lower_to_symmetric(rows: list[list[int]]) -> list[list[int]]:
    n = len(rows)
    fail_unless(all(len(row) == i + 1 for i, row in enumerate(rows)),
                "Malformed lower triangle.")
    out = [[0] * n for _ in range(n)]
    for i, row in enumerate(rows):
        for j, value in enumerate(row):
            fail_unless(type(value) is int, "Matrix coefficients must be integers.")
            out[i][j] = out[j][i] = value
    return out


def exact_matmul(A: list[list[int]], B: list[list[int]]) -> list[list[int]]:
    if np is not None:
        # dtype=object is essential: these are unbounded Python integers.
        return (np.array(A, dtype=object) @ np.array(B, dtype=object)).tolist()
    BT = list(zip(*B))
    return [[sum(x * y for x, y in zip(row, col)) for col in BT] for row in A]


def verify_positive_definite(R: list[list[int]], triangle: list[list[int]]) -> int:
    n = len(R)
    fail_unless(len(triangle) == n and
                all(len(row) == i + 1 for i, row in enumerate(triangle)),
                "Malformed triangular congruence certificate.")
    fail_unless(all(type(v) is int for row in triangle for v in row),
                "Congruence entries must be integers.")
    fail_unless(all(triangle[i][i] > 0 for i in range(n)),
                "The triangular congruence matrix must be invertible.")
    B = [row + [0] * (n - len(row)) for row in triangle]
    BT = [list(row) for row in zip(*B)]
    H = exact_matmul(exact_matmul(B, R), BT)
    fail_unless(all(H[i][j] == H[j][i] for i in range(n) for j in range(i)),
                "The congruence result must be symmetric.")
    margins = [H[i][i] - sum(abs(H[i][j]) for j in range(n) if j != i)
               for i in range(n)]
    fail_unless(min(margins) > 0,
                "Positive diagonal dominance failed.")
    # H is symmetric and strictly diagonally dominant with positive diagonal.
    # Therefore H is positive definite. B is invertible, so R is too.
    return min(margins)


def target_polynomial() -> Counter:
    """Formal contraction expansion of 4 K(T) + 3 V(T)^2 - 3 a^4.

    a = T_000, V(T) = ||T||_F^2 - a^2.
    Each nonzero dummy index ranges over the orthogonal complement of e_0.
    """
    target = Counter()
    tetrahedron = ((1, 2, 4), (1, 3, 5), (2, 3, 6), (4, 5, 6))
    for active in itertools.product((False, True), repeat=6):
        rows = tuple(tuple(label if active[label - 1] else 0 for label in row)
                     for row in tetrahedron)
        target[canonical_graph(rows)] += 4
    # Seven nonconstant tensor components; expand the square of their norm sum.
    for S in itertools.product((False, True), repeat=3):
        if not any(S):
            continue
        for T in itertools.product((False, True), repeat=3):
            if not any(T):
                continue
            first = tuple(i + 1 if S[i] else 0 for i in range(3))
            second = tuple(i + 4 if T[i] else 0 for i in range(3))
            target[canonical_graph((first, first, second, second))] += 3
    target[((0, 0, 0),) * 4] -= 3
    return target


def verify(path: Path) -> None:
    start = time.monotonic()
    certificate = json.loads(path.read_text())
    denominator = certificate["denominator"]
    fail_unless(type(denominator) is int and denominator > 0,
                "The denominator must be a positive integer.")
    actual = Counter()
    fail_unless([b["r"] for b in certificate["blocks"]] == list(range(7)),
                "The seven free-index blocks must be present.")
    for block in certificate["blocks"]:
        r = block["r"]
        features = tuple(tuple(tuple(row) for row in f) for f in block["features"])
        for feature in features:
            validate_feature(r, feature)
        N = block["combination_columns"]
        R = lower_to_symmetric(block["gram_lower"])
        fail_unless(len(N) == len(R), "Gram/combinations dimension mismatch.")
        for column in N:
            fail_unless(all(type(i) is int and type(a) is int and
                            0 <= i < len(features) for i, a in column),
                        "Invalid integer combination of quadratic features.")
        margin = verify_positive_definite(R, block["congruence_lower"])
        print(f"Block r={r}: exact positive definiteness passed ({len(R)} x {len(R)}).",
              flush=True)
        # Expand sum_{u,v} R_uv <sum_i N_iu F_i, sum_j N_jv F_j>.
        for u in range(len(R)):
            for v in range(u, len(R)):
                coefficient = R[u][v] * (1 if u == v else 2)
                if not coefficient:
                    continue
                for i, a in N[u]:
                    for j, b in N[v]:
                        graph = inner_product_graph(r, features[i], features[j])
                        actual[graph] += coefficient * a * b
    target = target_polynomial()
    keys = set(actual) | set(target)
    differences = {key: actual[key] - denominator * target[key]
                   for key in keys if actual[key] != denominator * target[key]}
    fail_unless(not differences,
                f"The exact polynomial identity failed at {len(differences)} coefficients.")
    print(f"All {len(keys)} formal contraction coefficients match exactly.")
    print("VERIFIED: 4 K(T) + 3 (||T||_F^2 - T_000^2)^2 - 3 T_000^4 >= 0")
    print("The identity and positivity certificates are independent of tensor dimension.")
    print(f"Elapsed: {time.monotonic() - start:.2f} seconds.")


def main() -> None:
    parser = argparse.ArgumentParser(description=__doc__)
    parser.add_argument("certificate", type=Path, nargs="?",
                        default=Path(__file__).with_name("certificate.json"))
    args = parser.parse_args()
    verify(args.certificate)


if __name__ == "__main__":
    main()
\end{lstlisting}

\end{document}